\documentclass[11pt]{amsart}
\newtheorem{theorem}{Theorem}[section]
\newtheorem{proposition}[theorem]{Proposition}
\newtheorem{corollary}[theorem]{Corollary}
\newtheorem{lemma}[theorem]{Lemma}
\theoremstyle{definition}
\newtheorem{definition}[theorem]{Definition}

\begin{document}

\title{Eigenvalue estimates for shrinkers}
\author{Simon Brendle and Raphael Tsiamis}
\address{Columbia University, 2990 Broadway, New York NY 10027, USA}
\address{Columbia University, 2990 Broadway, New York NY 10027, USA}
\thanks{The authors are grateful to Professor Gerhard Huisken for discussions. The first author was supported by the National Science Foundation under grant DMS-2103573 and by the Simons Foundation. He acknowledges the hospitality of T\"ubingen University, where part of this work was carried out. The second author was supported by an Onassis Foundation Scholarship and an A.G.~Leventis Foundation Scholarship.}

\begin{abstract}
We prove an eigenvalue estimate which holds on every properly embedded shrinker for mean curvature flow. This generalizes earlier work of Ding and Xin to the noncompact case.
\end{abstract}
\maketitle

\section{Introduction}

In this paper, we prove an eigenvalue estimate for properly embedded self-shrinkers of mean curvature flow. A hypersurface $M$ in $\mathbb{R}^{n+1}$ is called a self-shrinker if $H = \frac{1}{2} \, \langle x,\nu \rangle$, where $\nu$ denotes the unit normal vector field along $M$ and $H$ denotes the mean curvature of $M$. This condition implies that the hypersurface $M$ evolves under the mean curvature flow by shrinking homothetically. Shrinkers play a central role in understanding singularity formation for mean curvature flow (see, e.g., \cite{Colding-Minicozzi1}, \cite{Huisken}). 

\begin{theorem}
\label{main.thm}
Let $M$ be a properly embedded self-shrinker in $\mathbb{R}^{n+1}$, so that $H = \frac{1}{2} \, \langle x,\nu \rangle$. Suppose that $f \in H_{\text{\rm loc}}^1(M)$ satisfies $\int_M e^{-\frac{|x|^2}{4}} \, (f^2 + |\nabla^M f|^2) < \infty$ and $\int_M e^{-\frac{|x|^2}{4}} \, f = 0$. Then 
\[\int_M e^{-\frac{|x|^2}{4}} \, |\nabla^M f|^2 \geq \frac{1}{4} \int_M e^{-\frac{|x|^2}{4}} \, f^2.\] 
\end{theorem}

In the special case when $M$ is compact, Theorem \ref{main.thm} was proved by Ding and Xin (see \cite{Ding-Xin}, Theorem 1.3). Theorem \ref{main.thm} can be viewed as an analogue of the eigenvalue estimate for embedded minimal hypersurfaces proved by Choi and Wang \cite{Choi-Wang}. The eigenvalue estimate of Choi and Wang was later used by Choi and Schoen \cite{Choi-Schoen} to show that the space of compact embedded minimal surfaces of a bounded topology in $S^3$ is compact. We note that Colding and Minicozzi \cite{Colding-Minicozzi2} have obtained a compactness result for properly embedded shrinkers in $\mathbb{R}^3$ with bounded topology. 

The eigenvalue estimate in Theorem \ref{main.thm} is also related to the Frankel property for shrinkers; see \cite{Colding-Minicozzi3} and \cite{Impera-Pigola-Rimoldi} for further discussion.

The proof of Theorem \ref{main.thm} relies on a Reilly formula for the drift Laplacian. The drift Laplacian arises naturally in the context of shrinkers; in particular, it plays an important role in the fundamental work of Colding and Minicozzi \cite{Colding-Minicozzi1}. To handle the noncompact case, we minimize a certain functional, which depends on a parameter $\alpha>0$ (see Definition \ref{definition.of.mu} below). For each $\alpha>0$, we establish the existence and regularity of a minimizer (see Proposition \ref{existence.of.minimizer} and Proposition \ref{regularity}). Using a Reilly-type formula, we show that, for each $\alpha>0$, the minimum value of the functional is at least $\frac{1}{4}$ (see Proposition \ref{eigenvalue.estimate}). Theorem \ref{main.thm} follows by sending $\alpha \to 0$.

\section{A variational problem}

Throughout this paper, we fix a smooth function $\beta: [0,\infty) \to [0,1]$ such that $\beta(t) = 1$ for $t \in [0,1]$, $\beta'(t) \leq 0$ for $t \in [1,4]$, and $\beta(t) = 0$ for $t \in [4,\infty)$. For each positive integer $j$, we define a smooth function $\eta_j: \mathbb{R}^{n+1} \to [0,1]$ by 
\begin{equation} 
\label{definition.of.eta_j}
\eta_j(x) = \beta(j^{-2} \, |x|^2) 
\end{equation}
for all $x \in \mathbb{R}^{n+1}$.

\begin{definition}
Let $\mathcal{E}$ denote the set of all functions $w \in H_{\text{\rm loc}}^1(\mathbb{R}^{n+1})$ such that $\int_{\mathbb{R}^{n+1}} e^{-\frac{|x|^2}{4}} \, (w^2+|\nabla w|^2) < \infty$.  
\end{definition}

The following lemma is well known.

\begin{lemma}
\label{x.times.w.in.L2} 
Assume that $w \in H_{\text{\rm loc}}^1(\mathbb{R}^{n+1})$ satisfies $\int_{\mathbb{R}^{n+1}} e^{-\frac{|x|^2}{4}} \, |\nabla w|^2 < \infty$. Then 
\[\int_{\mathbb{R}^{n+1}} e^{-\frac{|x|^2}{4}} \, |x|^2 \, w^2 \leq C(n) \int_{B_{4n}(0)} w^2 + C(n) \int_{\mathbb{R}^{n+1}} e^{-\frac{|x|^2}{4}} \, |\nabla w|^2.\] 
In particular, $w \in \mathcal{E}$.
\end{lemma} 

\textbf{Proof.} 
Let $\eta_j$ denote the cutoff function defined in (\ref{definition.of.eta_j}). Using the inequality 
\[\langle x,\nabla \eta_j \rangle = 2 \, j^{-2} \, |x|^2 \, \beta'(j^{-2} \, |x|^2) \leq 0,\] 
we obtain 
\begin{align*} 
&\text{\rm div}(e^{-\frac{|x|^2}{4}} \, \eta_j^2 \, w^2 \, x) \\ 
&\leq \eta_j^2 \, \text{\rm div}(e^{-\frac{|x|^2}{4}} \, w^2 \, x) \\ 
&= e^{-\frac{|x|^2}{4}} \, \Big ( n+1 - \frac{|x|^2}{2} \Big ) \, \eta_j^2 \, w^2 + 2 \, e^{-\frac{|x|^2}{4}} \, \eta_j^2 \, w \, \langle x,\nabla w \rangle \\ 
&\leq e^{-\frac{|x|^2}{4}} \, \Big ( n+1 - \frac{|x|^2}{4} \Big ) \, \eta_j^2 \, w^2 + 4 \, e^{-\frac{|x|^2}{4}} \, \eta_j^2 \, |\nabla w|^2. 
\end{align*} 
This implies 
\begin{align*} 
0 &= \int_{\mathbb{R}^{n+1}} \text{\rm div}(e^{-\frac{|x|^2}{4}} \, \eta_j^2 \, w^2 \, x) \\ 
&\leq \int_{\mathbb{R}^{n+1}} e^{-\frac{|x|^2}{4}} \, \Big ( n+1 - \frac{|x|^2}{4} \Big ) \, \eta_j^2 \, w^2 + 4 \int_{\mathbb{R}^{n+1}} e^{-\frac{|x|^2}{4}} \, \eta_j^2 \, |\nabla w|^2. 
\end{align*}
Consequently, 
\[\int_{\mathbb{R}^{n+1}} e^{-\frac{|x|^2}{4}} \, |x|^2 \, \eta_j^2 \, w^2 \leq C(n) \int_{B_{4n}(0)} w^2 + C(n) \int_{\mathbb{R}^{n+1}} e^{-\frac{|x|^2}{4}} \, \eta_j^2 \, |\nabla w|^2.\] 
The claim follows by sending $j \to \infty$. This completes the proof of Lemma \ref{x.times.w.in.L2}. \\

In the following, we assume that $M$ is a properly embedded hypersurface in $\mathbb{R}^{n+1}$. The complement $\mathbb{R}^{n+1} \setminus M$ has two connected components, which we denote by $\Omega$ and $\tilde{\Omega}$. We denote by $\nu$ the outward-pointing unit normal to $\Omega$ and by $H$ the mean curvature of $\Omega$. We assume that $M$ satisfies the self-shrinker equation $H = \frac{1}{2} \, \langle x,\nu \rangle$.

\begin{proposition}
\label{area.growth}
The area of $M \cap B_r(0)$ is bounded from above by $O(r^n)$ as $r \to \infty$. 
\end{proposition} 

\textbf{Proof.} 
This estimate was proved by Ding and Xin (see \cite{Ding-Xin}, Theorem 2.2). Alternatively, the estimate can be deduced from Brakke's local area bounds for mean curvature flow (see \cite{Brakke} or \cite{Ecker}, Proposition 4.9). \\

\begin{corollary} 
\label{finite.weighted.area}
We have $\int_M e^{-\frac{|x|^2}{4}} < \infty$. 
\end{corollary}

\textbf{Proof.} 
This follows immediately from Proposition \ref{area.growth}. \\

\begin{definition}
Let $\mathcal{H}$ denote the set of all functions $f \in H_{\text{\rm loc}}^1(M)$ such that $\int_M e^{-\frac{|x|^2}{4}} \, (f^2+|\nabla^M f|^2) < \infty$. Let $\mathcal{H}_0 \subset \mathcal{H}$ denote the set of all functions $f \in \mathcal{H}$ satisfying $\int_M e^{-\frac{|x|^2}{4}} \, f^2 = 1$ and $\int_M e^{-\frac{|x|^2}{4}} \, f = 0$.
\end{definition}

In view of Corollary \ref{finite.weighted.area}, every constant function on $M$ belongs to the space $\mathcal{H}$. 

The following lemma is similar to Lemma B.1 in \cite{Bernstein-Wang}. 

\begin{lemma} 
\label{x.times.f.in.L2}
Assume that $f \in H_{\text{\rm loc}}^1(M)$ satisfies $\int_M e^{-\frac{|x|^2}{4}} \, |\nabla^M f|^2 < \infty$. Then 
\[\int_M e^{-\frac{|x|^2}{4}} \, |x|^2 \, f^2 \leq C(n) \int_{M \cap B_{4n}(0)} f^2 + C(n) \int_M e^{-\frac{|x|^2}{4}} \, |\nabla^M f|^2.\] 
In particular, $f \in \mathcal{H}$.
\end{lemma}

\textbf{Proof.}
Note that 
\[\text{\rm div}_M(x^{\text{\rm tan}}) = n - H \, \langle x,\nu \rangle = n - \frac{\langle x,\nu \rangle^2}{2}\] 
on $M$. Let $\eta_j$ denote the cutoff function defined in (\ref{definition.of.eta_j}). Using the inequality 
\[\langle x^{\text{\rm tan}},\nabla^M \eta_j \rangle = 2 \, j^{-2} \, |x^{\text{\rm tan}}|^2 \, \beta'(j^{-2} \, |x|^2) \leq 0,\] 
we obtain 
\begin{align*} 
&\text{\rm div}_M(e^{-\frac{|x|^2}{4}} \, \eta_j^2 \, f^2 \, x^{\text{\rm tan}}) \\ 
&\leq \eta_j^2 \, \text{\rm div}_M(e^{-\frac{|x|^2}{4}} \, f^2 \, x^{\text{\rm tan}}) \\ 
&= e^{-\frac{|x|^2}{4}} \, \Big ( n - \frac{\langle x,\nu \rangle^2}{2} - \frac{|x^{\text{\rm tan}}|^2}{2} \Big ) \, \eta_j^2 \, f^2 + 2 \, e^{-\frac{|x|^2}{4}} \, \eta_j^2 \, f \, \langle x^{\text{\rm tan}},\nabla^M f \rangle \\ 
&= e^{-\frac{|x|^2}{4}} \, \Big ( n - \frac{|x|^2}{2} \Big ) \, \eta_j^2 \, f^2 + 2 \, e^{-\frac{|x|^2}{4}} \, \eta_j^2 \, f \, \langle x^{\text{\rm tan}},\nabla^M f \rangle \\ 
&\leq e^{-\frac{|x|^2}{4}} \, \Big ( n - \frac{|x|^2}{4} \Big ) \, \eta_j^2 \, f^2 + 4 \, e^{-\frac{|x|^2}{4}} \, \eta_j^2 \, |\nabla^M f|^2 
\end{align*} 
on $M$. This implies
\begin{align*} 
0 &= \int_M \text{\rm div}_M(e^{-\frac{|x|^2}{4}} \, \eta_j^2 \, f^2 \, x^{\text{\rm tan}}) \\ 
&\leq \int_M e^{-\frac{|x|^2}{4}} \, \Big ( n - \frac{|x|^2}{4} \Big ) \, \eta_j^2 \, f^2 + 4 \int_M e^{-\frac{|x|^2}{4}} \, \eta_j^2 \, |\nabla^M f|^2. 
\end{align*} 
Consequently, 
\[\int_M e^{-\frac{|x|^2}{4}} \, |x|^2 \, \eta_j^2 \, f^2 \leq C(n) \int_{M \cap B_{4n}(0)} f^2 + C(n) \int_M e^{-\frac{|x|^2}{4}} \, \eta_j^2 \, |\nabla^M f|^2.\] 
The claim follows by sending $j \to \infty$. This completes the proof of Lemma \ref{x.times.f.in.L2}. \\

\begin{lemma}
\label{class.is.non.empty}
We can find a pair $(f,w) \in \mathcal{H}_0 \times \mathcal{E}$ such that $w|_M = f$.
\end{lemma}

\textbf{Proof.} 
We can find a function $w \in C_c^\infty(\mathbb{R}^{n+1})$ such that the restriction $w|_M$ is not constant. We define $w|_M = f$. Moreover, we put 
\[a = \frac{\int_M e^{-\frac{|x|^2}{4}} \, f}{\int_M e^{-\frac{|x|^2}{4}}}.\] 
Finally, we define $\hat{f} = f-a$ and $\hat{w} = w-a$. In view of Corollary \ref{finite.weighted.area}, constant functions on $M$ belong to $\mathcal{H}$. Therefore, $(\hat{f},\hat{w}) \in \mathcal{H} \times \mathcal{E}$ and $\hat{w}|_M = \hat{f}$. Moreover, $\int_M e^{-\frac{|x|^2}{4}} \, \hat{f}^2 > 0$ and $\int_M e^{-\frac{|x|^2}{4}} \, \hat{f} = 0$. Therefore, we may normalize the pair $(\hat{f},\hat{w}) \in \mathcal{H} \times \mathcal{E}$ by the constant $\int_M e^{-\frac{|x|^2}{4}} \hat{f}^2 > 0$ to obtain a pair in $\mathcal{H}_0 \times \mathcal{E}$ which has all the required properties. \\

\begin{definition}
\label{definition.of.mu}
Let us fix a real number $\alpha>0$. We define $\mu$ to be the infimum of the functional 
\[\int_M e^{-\frac{|x|^2}{4}} \, |\nabla^M f|^2 + \alpha \int_{\mathbb{R}^{n+1}} e^{-\frac{|x|^2}{4}} \, |\nabla w|^2\] 
over all pairs $(f,w) \in \mathcal{H}_0 \times \mathcal{E}$ satisfying $w|_M = f$. By Lemma \ref{class.is.non.empty}, the set of such pairs is non-empty. In particular, $\mu$ is finite.
\end{definition}

\begin{proposition}
\label{existence.of.minimizer}
Let us fix a real number $\alpha>0$. We can find a pair $(f,w) \in \mathcal{H}_0 \times \mathcal{E}$ such that $w|_M = f$ and 
\[\int_M e^{-\frac{|x|^2}{4}} \, |\nabla^M f|^2 + \alpha \int_{\mathbb{R}^{n+1}} e^{-\frac{|x|^2}{4}} \, |\nabla w|^2 = \mu.\] 
\end{proposition}

\textbf{Proof.} 
By definition of $\mu$, we can find a sequence $(f_k,w_k) \in \mathcal{H}_0 \times \mathcal{E}$ such that $w_k|_M = f_k$ for each $k$ and  
\begin{equation} 
\label{minimizing.sequence}
\int_M e^{-\frac{|x|^2}{4}} \, |\nabla^M f_k|^2 + \alpha \int_{\mathbb{R}^{n+1}} e^{-\frac{|x|^2}{4}} \, |\nabla w_k|^2 \to \mu 
\end{equation}
as $k \to \infty$.

Let us fix $r$ sufficiently large so that $r \geq 4n$ and $M \cap B_{\frac{r}{4}}(0) \neq \emptyset$. Since $f_k \in \mathcal{H}_0$ for each $k$, we obtain 
\begin{equation} 
\label{L2.bound.f}
\sup_k \int_{M \cap B_r(0)} f_k^2 < \infty. 
\end{equation} 
Combining (\ref{minimizing.sequence}), (\ref{L2.bound.f}), and Lemma \ref{x.times.f.in.L2}, we obtain 
\begin{equation} 
\label{L2.bound.for.x.times.f}
\sup_k \int_M e^{-\frac{|x|^2}{4}} \, |x|^2 \, f_k^2 < \infty. 
\end{equation}
Moreover, (\ref{minimizing.sequence}) implies 
\begin{equation} 
\label{L2.bound.for.gradient.w}
\sup_k \int_{B_r(0)} |\nabla w_k|^2 < \infty. 
\end{equation} 
By the Poincar\'e inequality, we can find a sequence of real numbers $c_k$ such that 
\begin{equation} 
\label{poincare}
\sup_k \int_{B_r(0)} (w_k-c_k)^2 < \infty. 
\end{equation}
Using (\ref{L2.bound.for.gradient.w}), (\ref{poincare}), and the Sobolev trace theorem, we deduce that 
\begin{equation} 
\label{sobolev.trace}
\sup_k \int_{M \cap B_{\frac{r}{2}}(0)} (f_k-c_k)^2 < \infty. 
\end{equation} 
Combining (\ref{L2.bound.f}) and (\ref{sobolev.trace}), we conclude that $\sup_k c_k^2 < \infty$. Using (\ref{poincare}), it follows that 
\begin{equation} 
\label{L2.bound.w}
\sup_k \int_{B_r(0)} w_k^2 < \infty. 
\end{equation} 
Combining (\ref{minimizing.sequence}), (\ref{L2.bound.w}), and Lemma \ref{x.times.w.in.L2}, we conclude that 
\begin{equation} 
\label{weighted.L2.bound.w}
\sup_k \int_{\mathbb{R}^{n+1}} e^{-\frac{|x|^2}{4}} \, w_k^2 < \infty. 
\end{equation}
After passing to a subsequence, we may assume that the sequence $f_k$ converges to a function $f$ weakly in $H_{\text{\rm loc}}^1(M)$, and the sequence $w_k$ converges to a function $w$ weakly in $H_{\text{\rm loc}}^1(\mathbb{R}^{n+1})$. Since $f_k \in \mathcal{H}_0$, we obtain 
\[\bigg | 1 - \int_M e^{-\frac{|x|^2}{4}} \, \eta_j^2 \, f_k^2 \bigg | = \bigg | \int_M e^{-\frac{|x|^2}{4}} \, (1-\eta_j^2) \, f_k^2 \bigg | \leq j^{-2} \int_M e^{-\frac{|x|^2}{4}} \, |x|^2 \, f_k^2\] 
and 
\[\bigg | \int_M e^{-\frac{|x|^2}{4}} \, \eta_j \, f_k \bigg | = \bigg | \int_M e^{-\frac{|x|^2}{4}} \, (1-\eta_j) \, f_k \bigg | \leq j^{-1} \int_M e^{-\frac{|x|^2}{4}} \, |x| \, |f_k|\] 
for all $j$ and $k$. In a first step, we fix $j$ and pass to the limit as $k \to \infty$; in a second step, we send $j \to \infty$. Using (\ref{L2.bound.for.x.times.f}), it follows that 
\[\int_M e^{-\frac{|x|^2}{4}} \, f^2 = 1\] 
and 
\[\int_M e^{-\frac{|x|^2}{4}} \, f = 0.\] 
Moreover, using (\ref{minimizing.sequence}), (\ref{weighted.L2.bound.w}), and standard lower semicontinuity arguments, we obtain   
\begin{align*} 
&\int_M e^{-\frac{|x|^2}{4}} \, |\nabla^M f|^2 < \infty, \\ 
&\int_{\mathbb{R}^{n+1}} e^{-\frac{|x|^2}{4}} \, (w^2+|\nabla w|^2) < \infty,
\end{align*} 
and 
\begin{equation} 
\label{limit}
\int_M e^{-\frac{|x|^2}{4}} \, |\nabla^M f|^2 + \alpha \int_{\mathbb{R}^{n+1}} e^{-\frac{|x|^2}{4}} \, |\nabla w|^2 \leq \mu. 
\end{equation}
Thus, we conclude that $(f,w) \in \mathcal{H}_0 \times \mathcal{E}$ and $w|_M = f$. In particular, equality holds in (\ref{limit}). This completes the proof of Proposition \ref{existence.of.minimizer}. \\

\begin{corollary} 
\label{weak.solution}
Let us fix a real number $\alpha>0$. Suppose that $(f,w) \in \mathcal{H}_0 \times \mathcal{E}$ is the minimizer constructed in Proposition \ref{existence.of.minimizer}. Then 
\[\int_M e^{-\frac{|x|^2}{4}} \, \langle \nabla^M f,\nabla^M \varphi \rangle + \alpha \int_{\mathbb{R}^{n+1}} e^{-\frac{|x|^2}{4}} \, \langle \nabla w,\nabla \psi \rangle = \mu \int_M e^{-\frac{|x|^2}{4}} \, f \, \varphi\]  
for all pairs $(\varphi,\psi) \in \mathcal{H} \times \mathcal{E}$ satisfying $\psi|_M = \varphi$. 
\end{corollary} 

\textbf{Proof.} 
Let us consider a pair $(\varphi,\psi) \in \mathcal{H} \times \mathcal{E}$ satisfying $\psi|_M = \varphi$. Let 
\[a = \frac{\int_M e^{-\frac{|x|^2}{4}} \, \varphi}{\int_M e^{-\frac{|x|^2}{4}}}.\] 
For each $s \in \mathbb{R}$, we define $\hat{f}_s = f + s(\varphi-a)$ and $\hat{w}_s = w + s(\psi-a)$. Then $(\hat{f}_s,\hat{w}_s) \in \mathcal{H} \times \mathcal{E}$, $\hat{w}_s|_M = \hat{f}_s$, and $\int_M e^{-\frac{|x|^2}{4}} \, \hat{f}_s = 0$ for each $s \in \mathbb{R}$. By definition of $\mu$, we obtain 
\[\int_M e^{-\frac{|x|^2}{4}} |\nabla^M \hat{f}_s|^2 + \alpha \int_{\mathbb{R}^{n+1}} e^{-\frac{|x|^2}{4}} \, |\nabla \hat{w}_s|^2 \geq \mu \int_M e^{-\frac{|x|^2}{4}} \, \hat{f}_s^2\] 
for each $s \in \mathbb{R}$. Moreover, equality holds for $s=0$. Taking the derivative at $s=0$ gives 
\[\int_M e^{-\frac{|x|^2}{4}} \, \langle \nabla^M f,\nabla^M \varphi \rangle + \alpha \int_{\mathbb{R}^{n+1}} e^{-\frac{|x|^2}{4}} \, \langle \nabla w,\nabla \psi \rangle = \mu \int_M e^{-\frac{|x|^2}{4}} \, f \, (\varphi-a).\] 
Since $\int_M e^{-\frac{|x|^2}{4}} \, f = 0$, the assertion follows. This completes the proof of Corollary \ref{weak.solution}. \\

\section{Regularity of the minimizer}

Throughout this section, we fix a real number $\alpha>0$. Let $(f,w) \in \mathcal{H}_0 \times \mathcal{E}$ denote the minimizer constructed in Proposition \ref{existence.of.minimizer}. In particular, $w|_M = f$. It follows from Corollary \ref{weak.solution} that 
\begin{equation} 
\label{pde.for.w}
\int_{\mathbb{R}^{n+1}} e^{-\frac{|x|^2}{4}} \, \langle \nabla w,\nabla \psi \rangle = 0 
\end{equation} 
for all $\psi \in \mathcal{E}$ satisfying $\psi|_M=0$. The standard regularity theory for elliptic PDE implies that $w \in C_{\text{\rm loc}}^\infty(\mathbb{R}^{n+1} \setminus M)$. 

As above, we denote by $\Omega$ and $\tilde{\Omega}$ the connected components of $\mathbb{R}^{n+1} \setminus M$. We denote by $\nu$ the outward-pointing unit normal to $\Omega$, by $h$ the second fundamental form of $\Omega$, and by $H$ the mean curvature of $\Omega$. Similarly, we denote by $\tilde{\nu} = -\nu$ the outward-pointing unit normal to $\tilde{\Omega}$, by $\tilde{h} = -h$ the second fundamental form of $\tilde{\Omega}$, and by $\tilde{H} = -H$ the mean curvature of $\tilde{\Omega}$. For brevity, we write $u=w|_\Omega$ and $\tilde{u}=w|_{\tilde{\Omega}}$. It follows from (\ref{pde.for.w}) that 
\[\Delta u - \frac{1}{2} \, \langle x,\nabla u \rangle = 0\] 
in $\Omega$ and 
\[\Delta \tilde{u} - \frac{1}{2} \, \langle x,\nabla \tilde{u} \rangle = 0\] 
in $\tilde{\Omega}$. We now state the main result of this section.

\begin{proposition}
\label{regularity}
The function $f$ is smooth. Moreover, $u$ is smooth up to the boundary of $\Omega$ and $\tilde{u}$ is smooth up to the boundary of $\tilde{\Omega}$. 
\end{proposition}

In the remainder of this section, we will give the proof of Proposition \ref{regularity}. To that end, we fix an arbitrary point $p \in M$. Moreover, we fix $\bar{r}$ sufficiently small so that $M \cap B_{\bar{r}}(p)$ is a graph over $T_p M$ with small slope. 

\begin{lemma}
\label{improved.regularity.f}
Let $s \geq 1$ be an integer. Suppose that $f \in H^s(M \cap B_r(p))$, $u \in H^s(\Omega \cap B_r(p))$, and $\tilde{u} \in H^s(\tilde{\Omega} \cap B_r(p))$ for each $r \in (0,\bar{r})$. Then $f \in H^{s+1}(M \cap B_r(p))$ for each $r \in (0,\bar{r})$.
\end{lemma}

\textbf{Proof.} 
Let us fix a real number $r_0 \in (0,\bar{r})$. We choose real numbers $r_1$ and $r_2$ such that $r_0 < r_1 < r_2 < \bar{r}$. We can find an open domain $D$ with smooth boundary such that $D \subset \Omega \cap B_{r_2}(p)$ and $\Omega \cap B_{r_1}(p) \subset D$. Similarly, we can find an open domain $\tilde{D}$ with smooth boundary such that $\tilde{D} \subset \tilde{\Omega} \cap B_{r_2}(p)$ and $\tilde{\Omega} \cap B_{r_1}(p) \subset \tilde{D}$. Then $\partial D \cap B_{r_1}(p) = \partial \tilde{D} \cap B_{r_1}(p) = M \cap B_{r_1}(p)$. Let $\zeta$ be a smooth radial cutoff function such that $\zeta$ is supported in $B_{r_1}(p)$ and $\zeta=1$ on $B_{r_0}(p)$. 

We define a function $g \in H^s(\partial D)$ by $g = \zeta f$ on $\partial D \cap B_{r_1}(p)$ and $g = 0$ on $\partial D \setminus B_{r_1}(p)$. Similarly, we define a function $\tilde{g} \in H^s(\partial \tilde{D})$ by $\tilde{g} = \zeta f$ on $\partial \tilde{D} \cap B_{r_1}(p)$ and $\tilde{g} = 0$ on $\partial \tilde{D} \setminus B_{r_1}(p)$. Let $v \in H^s(D)$ denote the harmonic extension of the function $g$ to $D$. Similarly, let $\tilde{v} \in H^s(\tilde{D})$ denote the harmonic extension of $\tilde{g}$ to $\tilde{D}$. We define $z \in H^s(D)$ by $z = \zeta u - v$, and we define $\tilde{z} \in H^s(\tilde{D})$ by $\tilde{z} = \zeta \tilde{u} - \tilde{v}$. 

Using the identity $\Delta u - \frac{1}{2} \, \langle x,\nabla u \rangle = 0$, we obtain 
\[\Delta z = \frac{1}{2} \, \zeta \, \langle x,\nabla u \rangle + 2 \, \langle \nabla \zeta,\nabla u \rangle + \Delta \zeta \, u\] 
in $D$. The expression on the right hand side belongs to the Sobolev space $H^{s-1}(D)$. Since the boundary trace of $z$ vanishes, it follows that $z \in H^{s+1}(D)$. Since $s \geq 1$, it follows that the normal derivative of $z$ is well-defined and belongs to the Sobolev space $H^{s-1}(\partial D)$. Moreover, 
\begin{equation} 
\label{integration.by.parts}
\int_{D \cap B_{r_0}(p)} \langle \nabla z,\nabla \psi \rangle + \int_{D \cap B_{r_0}(p)} \Delta z \, \psi = \int_{\partial D \cap B_{r_0}(p)} \langle \nabla z,\nu \rangle \, \psi 
\end{equation}
for every test function $\psi \in C_c^\infty(B_{r_0}(p))$. We next observe that $\Delta z = \frac{1}{2} \, \langle x,\nabla u \rangle$ on $D \cap B_{r_0}(p)$. Consequently, the identity (\ref{integration.by.parts}) can be rewritten as 
\begin{equation} 
\label{int.by.parts}
\int_{D \cap B_{r_0}(p)} \langle \nabla z,\nabla \psi \rangle + \frac{1}{2} \int_{D \cap B_{r_0}(p)} \langle x,\nabla u \rangle \, \psi = \int_{\partial D \cap B_{r_0}(p)} \langle \nabla z,\nu \rangle \, \psi 
\end{equation}
for every test function $\psi \in C_c^\infty(B_{r_0}(p))$. In the next step, we consider the Dirichlet-to-Neumann map on the domain $D$. It is well-known that the Dirichlet-to-Neumann map is a pseudo-differential operator of order $1$ (cf. \cite{Lee-Uhlmann}, Section 1); in particular, it defines a bounded linear operator $\mathcal{N}: H^s(\partial D) \to H^{s-1}(\partial D)$. Since $g \in H^s(\partial D)$, $\mathcal{N}(g)$ belongs to the Sobolev space $H^{s-1}(\partial D)$. Moreover, the function $\mathcal{N}(g)$ satisfies 
\begin{equation} 
\label{definition.of.D.to.N.map}
\int_{D \cap B_{r_0}(p)} \langle \nabla v,\nabla \psi \rangle = \int_{\partial D \cap B_{r_0}(p)} \mathcal{N}(g) \, \psi 
\end{equation}
for every test function $\psi \in C_c^\infty(B_{r_0}(p))$. Adding (\ref{int.by.parts}) and (\ref{definition.of.D.to.N.map}) gives 
\begin{align*} 
\int_{D \cap B_{r_0}(p)} \langle \nabla u,\nabla \psi \rangle + \frac{1}{2} \int_{D \cap B_{r_0}(p)} \langle x,\nabla u \rangle \, \psi 
&= \int_{\partial D \cap B_{r_0}(p)} \langle \nabla z,\nu \rangle \, \psi \\ 
&+ \int_{\partial D \cap B_{r_0}(p)} \mathcal{N}(g) \, \psi 
\end{align*}
for every test function $\psi \in C_c^\infty(B_{r_0}(p))$. Replacing $\psi$ by $e^{-\frac{|x|^2}{4}} \, \psi$ gives  
\begin{align} 
\label{identity.on.D}
\int_{D \cap B_{r_0}(p)} e^{-\frac{|x|^2}{4}} \, \langle \nabla u,\nabla \psi \rangle 
&= \int_{\partial D \cap B_{r_0}(p)} e^{-\frac{|x|^2}{4}} \, \langle \nabla z,\nu \rangle \, \psi \notag \\ 
&+ \int_{\partial D \cap B_{r_0}(p)} e^{-\frac{|x|^2}{4}} \, \mathcal{N}(g) \, \psi 
\end{align} 
for every test function $\psi \in C_c^\infty(B_{r_0}(p))$. Let $\tilde{\mathcal{N}}: H^s(\partial \tilde{D}) \to H^{s-1}(\partial \tilde{D})$ denote the Dirichlet-to-Neumann map on the domain $\tilde{D}$. An analogous argument gives 
\begin{align} 
\label{identity.on.tilde.D}
\int_{\tilde{D} \cap B_{r_0}(p)} e^{-\frac{|x|^2}{4}} \, \langle \nabla \tilde{u},\nabla \psi \rangle 
&= \int_{\partial \tilde{D} \cap B_{r_0}(p)} e^{-\frac{|x|^2}{4}} \, \langle \nabla \tilde{z},\tilde{\nu} \rangle \, \psi \notag \\ 
&+ \int_{\partial \tilde{D} \cap B_{r_0}(p)} e^{-\frac{|x|^2}{4}} \, \tilde{\mathcal{N}}(\tilde{g}) \, \psi 
\end{align}
for every test function $\psi \in C_c^\infty(B_{r_0}(p))$. Adding (\ref{identity.on.D}) and (\ref{identity.on.tilde.D}), we obtain 
\begin{align*} 
\int_{B_{r_0}(p)} e^{-\frac{|x|^2}{4}} \, \langle \nabla w,\nabla \psi \rangle 
&= \int_{M \cap B_{r_0}(p)} e^{-\frac{|x|^2}{4}} \, (\langle \nabla z,\nu \rangle + \langle \nabla \tilde{z},\tilde{\nu} \rangle) \, \psi \\ 
&+ \int_{M \cap B_{r_0}(p)} e^{-\frac{|x|^2}{4}} \, (\mathcal{N}(g) + \tilde{\mathcal{N}}(\tilde{g})) \, \psi 
\end{align*}
for every test function $\psi \in C_c^\infty(B_{r_0}(p))$. Using Corollary \ref{weak.solution}, we conclude that 
\begin{align*} 
&\int_M e^{-\frac{|x|^2}{4}} \, \langle \nabla^M f,\nabla^M \psi \rangle - \mu \int_M e^{-\frac{|x|^2}{4}} \, f \, \psi \\ 
&+ \alpha \int_{M \cap B_{r_0}(p)} e^{-\frac{|x|^2}{4}} \, (\langle \nabla z,\nu \rangle + \langle \nabla \tilde{z},\tilde{\nu} \rangle) \, \psi \\ 
&+ \alpha \int_{M \cap B_{r_0}(p)} e^{-\frac{|x|^2}{4}} \, (\mathcal{N}(g) + \tilde{\mathcal{N}}(\tilde{g})) \, \psi = 0 
\end{align*} 
for every test function $\psi \in C_c^\infty(B_{r_0}(p))$. In other words, the function $f$ satisfies 
\begin{align*} 
&\Delta_M f - \frac{1}{2} \, \langle x^{\text{\rm tan}},\nabla^M f \rangle + \mu \, f \\ 
&- \alpha \, (\langle \nabla z,\nu \rangle + \langle \nabla \tilde{z},\tilde{\nu} \rangle) - \alpha \, (\mathcal{N}(g) + \tilde{\mathcal{N}}(\tilde{g})) = 0 
\end{align*}
on $M \cap B_{r_0}(p)$, where the identity is understood in the sense of distributions. Note that $\langle \nabla z,\nu \rangle$ and $\mathcal{N}(g)$ belong to the Sobolev space $H^{s-1}(\partial D)$, and $\langle \nabla \tilde{z},\tilde{\nu} \rangle$ and $\tilde{\mathcal{N}}(\tilde{g})$ belong to the Sobolev space $H^{s-1}(\partial \tilde{D})$. Standard regularity results for elliptic PDE imply that $f \in H^{s+1}(M \cap B_r(p))$ for each $r \in (0,r_0)$. This completes the proof of Lemma \ref{improved.regularity.f}. \\

\begin{lemma}
\label{improved.regularity.u}
Let $s \geq 1$ be an integer. Suppose that $f \in H^{s+1}(M \cap B_r(p))$, $u \in H^s(\Omega \cap B_r(p))$, and $\tilde{u} \in H^s(\tilde{\Omega} \cap B_r(p))$ for each $r \in (0,\bar{r})$. Then $u \in H^{s+1}(\Omega \cap B_r(p))$ and $\tilde{u} \in H^{s+1}(\tilde{\Omega} \cap B_r(p))$ for each $r \in (0,\bar{r})$. 
\end{lemma}

\textbf{Proof.} 
Let us fix a real number $r_0 \in (0,\bar{r})$. As above, we choose real numbers $r_1$ and $r_2$ such that $r_0 < r_1 < r_2 < \bar{r}$. We can find an open domain $D$ with smooth boundary such that $D \subset \Omega \cap B_{r_2}(p)$ and $\Omega \cap B_{r_1}(p) \subset D$. Similarly, we can find an open domain $\tilde{D}$ with smooth boundary such that $\tilde{D} \subset \tilde{\Omega} \cap B_{r_2}(p)$ and $\tilde{\Omega} \cap B_{r_1}(p) \subset \tilde{D}$. Then $\partial D \cap B_{r_1}(p) = \partial \tilde{D} \cap B_{r_1}(p) = M \cap B_{r_1}(p)$. Let $\zeta$ be a smooth radial cutoff function such that $\zeta$ is supported in $B_{r_1}(p)$ and $\zeta=1$ on $B_{r_0}(p)$. 

We define a function $g \in H^{s+1}(\partial D)$ by $g = \zeta f$ on $\partial D \cap B_{r_1}(p)$ and $g = 0$ on $\partial D \setminus B_{r_1}(p)$. Similarly, we define a function $\tilde{g} \in H^{s+1}(\partial \tilde{D})$ by $\tilde{g} = \zeta f$ on $\partial \tilde{D} \cap B_{r_1}(p)$ and $\tilde{g} = 0$ on $\partial \tilde{D} \setminus B_{r_1}(p)$. Let $v \in H^{s+1}(D)$ denote the harmonic extension of the function $g$ to $D$. Similarly, let $\tilde{v} \in H^{s+1}(\tilde{D})$ denote the harmonic extension of $\tilde{g}$ to $\tilde{D}$. We define $z \in H^s(D)$ by $z = \zeta u - v$, and we define $\tilde{z} \in H^s(\tilde{D})$ by $\tilde{z} = \zeta \tilde{u} - \tilde{v}$. 

Arguing as above, we can show that $\Delta z$ belongs to the Sobolev space $H^{s-1}(D)$. Since the boundary trace of $z$ vanishes, it follows that $z \in H^{s+1}(D)$. Since $v \in H^{s+1}(D)$, it follows that $\zeta u \in H^{s+1}(D)$ and $u \in H^{s+1}(\Omega \cap B_{r_0}(p))$. An analogous argument shows that $\tilde{u} \in H^{s+1}(\tilde{\Omega} \cap B_{r_0}(p))$. This completes the proof of Lemma \ref{improved.regularity.u}. \\

Combining Lemma \ref{improved.regularity.f} and Lemma \ref{improved.regularity.u}, we conclude that $f \in H^s(M \cap B_r(p))$, $u \in H^s(\Omega \cap B_r(p))$, and $\tilde{u} \in H^s(\tilde{\Omega} \cap B_r(p))$ for each $r \in (0,\bar{r})$ and every integer $s \geq 1$. This shows that $f$ is smooth, $u$ is smooth up to the boundary of $\Omega$, and $\tilde{u}$ is smooth up to the boundary of $\tilde{\Omega}$.  This completes the proof of Proposition \ref{regularity}.

\section{A lower bound for the eigenvalue $\mu$}

In this section, we establish a lower bound for the eigenvalue $\mu$ introduced in Definition \ref{definition.of.mu}. As in the work of Ding and Xin \cite{Ding-Xin}, we use a Reilly formula. Our arguments are different from those in \cite{Ding-Xin} due to the presence of the parameter $\alpha$ and the use of cutoff functions.

\begin{proposition}
\label{eigenvalue.estimate}
Let us fix a real number $\alpha>0$, and let $\mu$ be defined as in Definition \ref{definition.of.mu}. Then $\mu \geq \frac{1}{4}$.
\end{proposition}

\textbf{Proof.} 
Let $(f,w) \in \mathcal{H}_0 \times \mathcal{E}$ denote the minimizer constructed in Proposition \ref{existence.of.minimizer}. In particular, $w|_M = f$. As above, we denote by $\Omega$ and $\tilde{\Omega}$ the connected components of $\mathbb{R}^{n+1} \setminus M$. Moreover, we write $u=w|_\Omega$ and $\tilde{u}=w|_{\tilde{\Omega}}$. By Proposition \ref{regularity}, $f$ is smooth, $u$ is smooth up to the boundary of $\Omega$, and $\tilde{u}$ is smooth up to the boundary of $\tilde{\Omega}$. Moreover, $\Delta u - \frac{1}{2} \, \langle x,\nabla u \rangle = 0$ in $\Omega$ and $\Delta \tilde{u} - \frac{1}{2} \, \langle x,\nabla \tilde{u} \rangle = 0$ in $\tilde{\Omega}$. Corollary \ref{weak.solution} implies that 
\begin{equation} 
\label{pde.for.f}
\Delta_M f - \frac{1}{2} \, \langle x^{\text{\rm tan}},\nabla^M f \rangle + \mu f - \alpha \, (\langle \nabla u,\nu \rangle + \langle \nabla \tilde{u},\tilde{\nu} \rangle) = 0 
\end{equation}
on $M$. 

Let $\eta_j$ denote the cutoff function defined in (\ref{definition.of.eta_j}). Using the identity $\Delta u - \frac{1}{2} \, \langle x,\nabla u \rangle = 0$, we obtain 
 \begin{align*} 
&\frac{1}{2} \, \text{\rm div} \big ( e^{-\frac{|x|^2}{4}} \, \nabla(|\nabla u|^2) \big ) \\ 
&= \frac{1}{2} \, e^{-\frac{|x|^2}{4}} \, \Delta(|\nabla u|^2) - \frac{1}{4} \, e^{-\frac{|x|^2}{4}} \, \big \langle x,\nabla(|\nabla u|^2) \big \rangle \\ 
&= e^{-\frac{|x|^2}{4}} \, |D^2 u|^2 + e^{-\frac{|x|^2}{4}} \, \langle \nabla u,\nabla \Delta u \rangle - \frac{1}{4} \, e^{-\frac{|x|^2}{4}} \, \big \langle x,\nabla(|\nabla u|^2) \big \rangle \\ 
&= e^{-\frac{|x|^2}{4}} \, |D^2 u|^2 + \frac{1}{2} \, e^{-\frac{|x|^2}{4}} \, \big \langle \nabla u,\nabla(\langle x,\nabla u \rangle) \big \rangle - \frac{1}{4} \, e^{-\frac{|x|^2}{4}} \, \big \langle x,\nabla(|\nabla u|^2) \big \rangle \\ 
&= e^{-\frac{|x|^2}{4}} \, |D^2 u|^2 + \frac{1}{2} \, e^{-\frac{|x|^2}{4}} \, |\nabla u|^2. 
\end{align*} 
We multiply this identity by $\eta_j^2$ and integrate over $\Omega$. We arrive at  
\begin{align} 
\label{divergence.thm.on.Omega}
&\frac{1}{2} \int_M e^{-\frac{|x|^2}{4}} \, \eta_j^2 \, \big \langle \nabla(|\nabla u|^2),\nu \big \rangle \notag \\ 
&= \int_\Omega e^{-\frac{|x|^2}{4}} \, \eta_j^2 \, |D^2 u|^2 + \frac{1}{2} \int_\Omega e^{-\frac{|x|^2}{4}} \, \eta_j^2 \, |\nabla u|^2 \\ 
&+ \int_\Omega e^{-\frac{|x|^2}{4}} \, \eta_j \, \big \langle \nabla \eta_j,\nabla(|\nabla u|^2) \big \rangle. \notag
\end{align}
On the other hand, using the identities $\Delta u - \frac{1}{2} \, \langle x,\nabla u \rangle = 0$ and $H = \frac{1}{2} \, \langle x,\nu \rangle$, we obtain 
\begin{align*} 
(D^2 u)(\nu,\nu) 
&= -\Delta_M f + \Delta u - H \, \langle \nabla u,\nu \rangle \\ 
&= -\Delta_M f + \frac{1}{2} \, \langle x,\nabla u \rangle - \frac{1}{2} \, \langle x,\nu \rangle \, \langle \nabla u,\nu \rangle \\ 
&= -\Delta_M f + \frac{1}{2} \, \langle x^{\text{\rm tan}},\nabla^M f \rangle 
\end{align*} 
at each point on $M$. This implies 
\begin{align*} 
&\frac{1}{2} \, e^{-\frac{|x|^2}{4}} \, \big \langle \nabla(|\nabla u|^2),\nu \rangle \\ 
&= e^{-\frac{|x|^2}{4}} \, (D^2 u)(\nabla u,\nu) \\ 
&= e^{-\frac{|x|^2}{4}} \, (D^2 u)(\nu,\nu) \, \langle \nabla u,\nu \rangle + e^{-\frac{|x|^2}{4}} \, (D^2 u)(\nabla^M f,\nu) \\ 
&= -e^{-\frac{|x|^2}{4}} \ \Big ( \Delta_M f - \frac{1}{2} \, \langle x^{\text{\rm tan}},\nabla^M f \rangle \Big ) \, \langle \nabla u,\nu \rangle \\ 
&+ e^{-\frac{|x|^2}{4}} \, \big \langle \nabla^M f,\nabla^M(\langle \nabla u,\nu \rangle) \big \rangle - e^{-\frac{|x|^2}{4}} \, h(\nabla^M f,\nabla^M f) \\ 
&= -2 \, e^{-\frac{|x|^2}{4}} \ \Big ( \Delta_M f - \frac{1}{2} \, \langle x^{\text{\rm tan}},\nabla^M f \rangle \Big ) \, \langle \nabla u,\nu \rangle \\ 
&+ \text{\rm div}_M \big ( e^{-\frac{|x|^2}{4}} \, \langle \nabla u,\nu \rangle \, \nabla^M f \big ) - e^{-\frac{|x|^2}{4}} \, h(\nabla^M f,\nabla^M f) 
\end{align*}
at each point on $M$. We multiply this identity by $\eta_j^2$ and integrate over $M$. This gives 
\begin{align} 
\label{divergence.thm.on.M}
&\frac{1}{2} \int_M e^{-\frac{|x|^2}{4}} \, \eta_j^2 \, \big \langle \nabla(|\nabla u|^2),\nu \big \rangle \notag \\ 
&= -2 \int_M e^{-\frac{|x|^2}{4}} \, \eta_j^2 \, \Big ( \Delta_M f - \frac{1}{2} \, \langle x^{\text{\rm tan}},\nabla^M f \rangle \Big ) \, \langle \nabla u,\nu \rangle \\ 
&- 2 \int_M e^{-\frac{|x|^2}{4}} \, \eta_j \, \langle \nabla^M \eta_j,\nabla^M f \rangle \, \langle \nabla u,\nu \rangle - \int_M e^{-\frac{|x|^2}{4}} \, \eta_j^2 \, h(\nabla^M f,\nabla^M f). \notag
\end{align}
Combining (\ref{divergence.thm.on.Omega}) and (\ref{divergence.thm.on.M}), we conclude that 
\begin{align} 
\label{formula.on.Omega}
&- 2 \int_M e^{-\frac{|x|^2}{4}} \, \eta_j^2 \, \Big ( \Delta_M f - \frac{1}{2} \, \langle x^{\text{\rm tan}},\nabla^M f \rangle \Big ) \, \langle \nabla u,\nu \rangle \notag \\ 
&- 2 \int_M e^{-\frac{|x|^2}{4}} \, \eta_j \, \langle \nabla^M \eta_j,\nabla^M f \rangle \, \langle \nabla u,\nu \rangle - \int_M e^{-\frac{|x|^2}{4}} \, \eta_j^2 \, h(\nabla^M f,\nabla^M f) \notag \\ 
&= \int_\Omega e^{-\frac{|x|^2}{4}} \, \eta_j^2 \, |D^2 u|^2 + \frac{1}{2} \int_\Omega e^{-\frac{|x|^2}{4}} \, \eta_j^2 \, |\nabla u|^2 \\ 
&+ \int_\Omega e^{-\frac{|x|^2}{4}} \, \eta_j \, \big \langle \nabla \eta_j,\nabla(|\nabla u|^2) \big \rangle. \notag
\end{align}
An analogous computation gives 
\begin{align} 
\label{formula.on.tilde.Omega}
&-2 \int_M e^{-\frac{|x|^2}{4}} \, \eta_j^2 \, \Big ( \Delta_M f + \frac{1}{2} \, \langle x^{\text{\rm tan}},\nabla^M f \rangle \Big ) \, \langle \nabla \tilde{u},\tilde{\nu} \rangle \notag \\ 
&- 2 \int_M e^{-\frac{|x|^2}{4}} \, \eta_j \, \langle \nabla^M \eta_j,\nabla^M f \rangle \, \langle \nabla \tilde{u},\tilde{\nu} \rangle - \int_M e^{-\frac{|x|^2}{4}} \, \eta_j^2 \, \tilde{h}(\nabla^M f,\nabla^M f) \notag \\ 
&= \int_{\tilde{\Omega}} e^{-\frac{|x|^2}{4}} \, \eta_j^2 \, |D^2 \tilde{u}|^2 + \frac{1}{2} \int_{\tilde{\Omega}} e^{-\frac{|x|^2}{4}} \, \eta_j^2 \, |\nabla \tilde{u}|^2 \\ 
&+ \int_{\tilde{\Omega}} e^{-\frac{|x|^2}{4}} \, \eta_j \, \big \langle \nabla \eta_j,\nabla(|\nabla \tilde{u}|^2) \big \rangle. \notag
\end{align}
In the next step, we add (\ref{formula.on.Omega}) and (\ref{formula.on.tilde.Omega}). Using the identity $h+\tilde{h}=0$, we obtain 
\begin{align} 
\label{sum}
&- 2 \int_M e^{-\frac{|x|^2}{4}} \, \eta_j^2 \, \Big ( \Delta_M f - \frac{1}{2} \, \langle x^{\text{\rm tan}},\nabla^M f \rangle \Big ) \, (\langle \nabla u,\nu \rangle + \langle \nabla \tilde{u},\tilde{\nu} \rangle) \notag \\ 
&- 2 \int_M e^{-\frac{|x|^2}{4}} \, \eta_j \, \langle \nabla^M \eta_j,\nabla^M f \rangle \, (\langle \nabla u,\nu \rangle + \langle \nabla \tilde{u},\tilde{\nu} \rangle) \notag \\ 
&= \int_\Omega e^{-\frac{|x|^2}{4}} \, \eta_j^2 \, |D^2 u|^2 + \int_{\tilde{\Omega}} e^{-\frac{|x|^2}{4}} \, \eta_j^2 \, |D^2 \tilde{u}|^2 \\ 
&+ \frac{1}{2} \int_\Omega e^{-\frac{|x|^2}{4}} \, \eta_j^2 \, |\nabla u|^2 + \frac{1}{2} \int_{\tilde{\Omega}} e^{-\frac{|x|^2}{4}} \, \eta_j^2 \, |\nabla \tilde{u}|^2 \notag \\ 
&+ \int_\Omega e^{-\frac{|x|^2}{4}} \, \eta_j \, \big \langle \nabla \eta_j,\nabla(|\nabla u|^2) \big \rangle + \int_{\tilde{\Omega}} e^{-\frac{|x|^2}{4}} \, \eta_j \, \big \langle \nabla \eta_j,\nabla(|\nabla \tilde{u}|^2) \big \rangle. \notag
\end{align}
Substituting (\ref{pde.for.f}) into (\ref{sum}) gives  
\begin{align} 
\label{key.identity}
&2\mu \int_M e^{-\frac{|x|^2}{4}} \, \eta_j^2 \, f \, (\langle \nabla u,\nu \rangle + \langle \nabla \tilde{u},\tilde{\nu} \rangle) \notag \\ 
&- 2\alpha \int_M e^{-\frac{|x|^2}{4}} \, \eta_j^2 \, (\langle \nabla u,\nu \rangle + \langle \nabla \tilde{u},\tilde{\nu} \rangle)^2 \notag \\ 
&- 2 \int_M e^{-\frac{|x|^2}{4}} \, \eta_j \, \langle \nabla^M \eta_j,\nabla^M f \rangle \, (\langle \nabla u,\nu \rangle + \langle \nabla \tilde{u},\tilde{\nu} \rangle) \notag \\ 
&= \int_\Omega e^{-\frac{|x|^2}{4}} \, \eta_j^2 \, |D^2 u|^2 + \int_{\tilde{\Omega}} e^{-\frac{|x|^2}{4}} \, \eta_j^2 \, |D^2 \tilde{u}|^2 \\ 
&+ \frac{1}{2} \int_\Omega e^{-\frac{|x|^2}{4}} \, \eta_j^2 \, |\nabla u|^2 + \frac{1}{2} \int_{\tilde{\Omega}} e^{-\frac{|x|^2}{4}} \, \eta_j^2 \, |\nabla \tilde{u}|^2 \notag \\ 
&+ \int_\Omega e^{-\frac{|x|^2}{4}} \, \eta_j \, \big \langle \nabla \eta_j,\nabla(|\nabla u|^2) \big \rangle + \int_{\tilde{\Omega}} e^{-\frac{|x|^2}{4}} \, \eta_j \, \big \langle \nabla \eta_j,\nabla(|\nabla \tilde{u}|^2) \big \rangle. \notag 
\end{align}
Using the divergence theorem, we obtain 
\begin{align} 
\label{div}
&\int_M e^{-\frac{|x|^2}{4}} \, \eta_j^2 \, f \, (\langle \nabla u,\nu \rangle + \langle \nabla \tilde{u},\tilde{\nu} \rangle) \notag \\ 
&= \int_\Omega e^{-\frac{|x|^2}{4}} \, \eta_j^2 \, |\nabla u|^2 + \int_{\tilde{\Omega}} e^{-\frac{|x|^2}{4}} \, \eta_j^2 \, |\nabla \tilde{u}|^2 \\ 
&+ 2 \int_\Omega e^{-\frac{|x|^2}{4}} \, \eta_j \, u \, \langle \nabla \eta_j,\nabla u \rangle + 2 \int_{\tilde{\Omega}} e^{-\frac{|x|^2}{4}} \, \eta_j \, \tilde{u} \, \langle \nabla \eta_j,\nabla \tilde{u} \rangle. \notag
\end{align} 
Substituting (\ref{div}) into (\ref{key.identity}), we conclude that 
\begin{align} 
\label{key.identity.2}
&\Big ( 2\mu - \frac{1}{2} \Big ) \int_\Omega e^{-\frac{|x|^2}{4}} \, \eta_j^2 \, |\nabla u|^2 + \Big ( 2\mu - \frac{1}{2} \Big ) \int_{\tilde{\Omega}} e^{-\frac{|x|^2}{4}} \, \eta_j^2 \, |\nabla \tilde{u}|^2 \notag \\ 
&= 2\alpha \int_M e^{-\frac{|x|^2}{4}} \, \eta_j^2 \, (\langle \nabla u,\nu \rangle + \langle \nabla \tilde{u},\tilde{\nu} \rangle)^2 \notag \\ 
&+ 2 \int_M e^{-\frac{|x|^2}{4}} \, \eta_j \, \langle \nabla^M \eta_j,\nabla^M f \rangle \, (\langle \nabla u,\nu \rangle + \langle \nabla \tilde{u},\tilde{\nu} \rangle) \notag \\ 
&+ \int_\Omega e^{-\frac{|x|^2}{4}} \, \eta_j^2 \, |D^2 u|^2 + \int_{\tilde{\Omega}} e^{-\frac{|x|^2}{4}} \, \eta_j^2 \, |D^2 \tilde{u}|^2 \\ 
&+ \int_\Omega e^{-\frac{|x|^2}{4}} \, \eta_j \, \big \langle \nabla \eta_j,\nabla(|\nabla u|^2) \big \rangle + \int_{\tilde{\Omega}} e^{-\frac{|x|^2}{4}} \, \eta_j \, \big \langle \nabla \eta_j,\nabla(|\nabla \tilde{u}|^2) \big \rangle \notag \\ 
&- 4\mu \int_\Omega e^{-\frac{|x|^2}{4}} \, \eta_j \, u \, \langle \nabla \eta_j,\nabla u \rangle - 4\mu \int_{\tilde{\Omega}} e^{-\frac{|x|^2}{4}} \, \eta_j \, \tilde{u} \, \langle \nabla \eta_j,\nabla \tilde{u} \rangle.  \notag
\end{align}
Using Young's inequality, we obtain 
\begin{align} 
\label{aux.1}
&2\alpha \, \eta_j^2 \, (\langle \nabla u,\nu \rangle + \langle \nabla \tilde{u},\tilde{\nu} \rangle)^2 + 2\eta_j \, \langle \nabla^M \eta_j,\nabla^M f \rangle \, (\langle \nabla u,\nu \rangle + \langle \nabla \tilde{u},\tilde{\nu} \rangle) \notag \\ 
&\geq -\frac{1}{2\alpha} \, \langle \nabla^M \eta_j,\nabla^M f \rangle^2 \\
&\geq -\frac{1}{2\alpha} \, |\nabla^M \eta_j|^2 \, |\nabla^M f|^2  \notag
\end{align}
at each point on $M$. Moreover, 
\begin{align}
\label{aux.2}
&\eta_j^2 \, |D^2 u|^2 + \eta_j \, \big \langle \nabla \eta_j,\nabla(|\nabla u|^2) \big \rangle \notag \\ 
&\geq \eta_j^2 \, |D^2 u|^2 - 2\eta_j \, |\nabla \eta_j| \, |\nabla u| \, |D^2 u| \\ 
&\geq -|\nabla \eta_j|^2 \, |\nabla u|^2 \notag  
\end{align}
at each point in $\Omega$ and 
\begin{align}\label{aux.3} 
&\eta_j^2 \, |D^2 \tilde{u}|^2 + \eta_j \, \big \langle \nabla \eta_j,\nabla(|\nabla \tilde{u}|^2) \big \rangle \notag \\ 
&\geq \eta_j^2 \, |D^2 \tilde{u}|^2 - 2\eta_j \, |\nabla \eta_j| \, |\nabla \tilde{u}| \, |D^2 \tilde{u}| \\
& \geq -|\nabla \eta_j|^2 \, |\nabla \tilde{u}|^2         \notag
    \end{align}
at each point in $\tilde{\Omega}$. Combining (\ref{key.identity.2}), (\ref{aux.1}), (\ref{aux.2}), and (\ref{aux.3}), we obtain 
\begin{align*} 
&\Big ( 2\mu - \frac{1}{2} \Big ) \int_\Omega e^{-\frac{|x|^2}{4}} \, \eta_j^2 \, |\nabla u|^2 + \Big ( 2\mu - \frac{1}{2} \Big ) \int_{\tilde{\Omega}} e^{-\frac{|x|^2}{4}} \, \eta_j^2 \, |\nabla \tilde{u}|^2 \\ 
&\geq -\frac{1}{2\alpha} \int_M e^{-\frac{|x|^2}{4}} \, |\nabla^M \eta_j|^2 \, |\nabla^M f|^2 \\ 
&- \int_\Omega e^{-\frac{|x|^2}{4}} \, |\nabla \eta_j|^2 \, |\nabla u|^2 - \int_{\tilde{\Omega}} e^{-\frac{|x|^2}{4}} \, |\nabla \eta_j|^2 \, |\nabla \tilde{u}|^2 \\ 
&- 4\mu \int_\Omega e^{-\frac{|x|^2}{4}} \, \eta_j \, u \, \langle \nabla \eta_j,\nabla u \rangle - 4\mu \int_{\tilde{\Omega}} e^{-\frac{|x|^2}{4}} \, \eta_j \, \tilde{u} \, \langle \nabla \eta_j,\nabla \tilde{u} \rangle. 
\end{align*}
Since $\int_M e^{-\frac{|x|^2}{4}} \, |\nabla^M f|^2 < \infty$ and $\int_{\mathbb{R}^{n+1}} e^{-\frac{|x|^2}{4}} \, (w^2 + |\nabla w|^2) < \infty$, the terms on the right hand side converge to $0$ as $j \to \infty$. We conclude that 
\[\Big ( 2\mu - \frac{1}{2} \Big ) \int_\Omega e^{-\frac{|x|^2}{4}} \, |\nabla u|^2 + \Big ( 2\mu - \frac{1}{2} \Big ) \int_{\tilde{\Omega}} e^{-\frac{|x|^2}{4}} \, |\nabla \tilde{u}|^2 \geq 0.\] 
Suppose that $\mu<\frac{1}{4}$. Then $\nabla u$ and $\nabla \tilde{u}$ vanish identically. This implies that $w$ is constant. Consequently, $f$ is constant. Since $\int_M e^{-\frac{|x|^2}{4}} \, f = 0$, it follows that $f$ vanishes identically. This contradicts the fact that $\int_M e^{-\frac{|x|^2}{4}} \, f^2 = 1$. This completes the proof of Proposition \ref{eigenvalue.estimate}. \\

\begin{corollary} 
\label{eigenvalue.estimate.2}
Assume that $f \in \mathcal{H}_0$. Moreover, suppose that there exists a function $w \in \mathcal{E}$ such that $w|_M = f$. Then 
\[\int_M e^{-\frac{|x|^2}{4}} \, |\nabla^M f|^2 \geq \frac{1}{4}.\] 
\end{corollary}

\textbf{Proof.} 
Proposition \ref{eigenvalue.estimate} implies that 
\[\int_M e^{-\frac{|x|^2}{4}} \, |\nabla^M f|^2 + \alpha \int_{\mathbb{R}^{n+1}} e^{-\frac{|x|^2}{4}} \, |\nabla w|^2 \geq \frac{1}{4}\] 
for each $\alpha>0$. Since $\alpha>0$ is arbitrary, the assertion follows. This completes the proof of Corollary \ref{eigenvalue.estimate.2}. \\

\section{Proof of Theorem \ref{main.thm}}

We now give the proof of Theorem \ref{main.thm}. Assume that $f \in \mathcal{H}$ satisfies $\int_M e^{-\frac{|x|^2}{4}} \, f = 0$. Let $\eta_j$ denote the cutoff function defined in (\ref{definition.of.eta_j}). For each $j$, we can find a function $w_j \in H_{\text{\rm loc}}^1(\mathbb{R}^{n+1})$ such that $w_j$ has compact support and $w_j|_M = \eta_j f$. We define a sequence of real numbers $a_j$ by 
\[a_j = \frac{\int_M e^{-\frac{|x|^2}{4}} \, \eta_j f}{\int_M e^{-\frac{|x|^2}{4}}}.\] 
Since $\int_M e^{-\frac{|x|^2}{4}} \, f = 0$, it follows that $a_j \to 0$. 

For each $j$, we define $\hat{f}_j = \eta_j f - a_j$ and $\hat{w}_j = w_j - a_j$. Clearly, $(\hat{f}_j,\hat{w}_j) \in \mathcal{H} \times \mathcal{E}$, $\hat{w}_j|_M = \hat{f}_j$, and $\int_M e^{-\frac{|x|^2}{4}} \, \hat{f}_j = 0$ for each $j$. Using Corollary \ref{eigenvalue.estimate.2}, we obtain 
\[\int_M e^{-\frac{|x|^2}{4}} \, |\nabla^M \hat{f}_j|^2 \geq \frac{1}{4} \int_M e^{-\frac{|x|^2}{4}} \, \hat{f}_j^2\] 
for each $j$. This implies 
\[\int_M e^{-\frac{|x|^2}{4}} \, |\nabla^M (\eta_j f)|^2 \geq \frac{1}{4} \int_M e^{-\frac{|x|^2}{4}} \, (\eta_j f - a_j)^2\] 
for each $j$. Sending $j \to \infty$ gives 
\[\int_M e^{-\frac{|x|^2}{4}} \, |\nabla^M f|^2 \geq \frac{1}{4} \int_M e^{-\frac{|x|^2}{4}} \, f^2.\] 
This completes the proof of Theorem \ref{main.thm}.

\end{document}